\documentclass{amsart}
\usepackage{pifont,color}
\IfFileExists{times.sty}{\usepackage{times}}{}
\newcounter{commentlabel}
\newcommand{\COMMENT}[1]{\stepcounter{commentlabel}\ifnum\thecommentlabel=212\setcounter{commentlabel}{171}\fi{\normalsize\textcolor{red}{\ding{\thecommentlabel}}}\marginpar{\color{red}\footnotesize\vskip-.6\baselineskip\noindent\raggedright\hsize1.1in\ding{\thecommentlabel}\thinspace#1\vskip.3\baselineskip}}
\newcommand{\SecDef}{2}      % Definitions
\newcommand{\SecPiece}{8}    % Piecewise linear maps
\theoremstyle{plain}
\newtheorem{theorem}{Theorem}
\newtheorem{lemma}[theorem]{Lemma}
\newtheorem{proposition}[theorem]{Proposition}
\newtheorem{corollary}[theorem]{Corollary}

\newtheorem{conjecture}[theorem]{Conjecture}
\theoremstyle{definition}
\newtheorem{definition}[theorem]{Definition}
\newtheorem{example}[theorem]{Example}
\newtheorem{remark}[theorem]{Remark}
\newtheorem{question}[theorem]{Question}
\numberwithin{equation}{section}
\numberwithin{theorem}{section}

\begin{document}

% Bibliographic entries
\newcommand{\AABM}{{AABM}}     % Abarenkova, Angles d'Auriac, Bourkraa 
                               % and Maillard
\newcommand{\AKM}{{AKM}}       % Adler, Konheim and McAndrew
\newcommand{\Ar}{{Ar}}         % Arnold
\newcommand{\BK}{{BK}}         % Bedford and Kim
\newcommand{\Be}{{Be}}         % Bellon
\newcommand{\BFJ}{{BFJ}}       % Boucksom, Favre, and Jonsson
\newcommand{\Bo}{{Bo}}         % Bowen
\newcommand{\BV}{{BV}}         % Bellon and Viallet
\newcommand{\BM}{{BM}}         % Bourkraa and Maillard
\newcommand{\DF}{{DF}}         % Diller and Favre
\newcommand{\Di}{{Di}}         % Dinaburg
\newcommand{\DS}{{DS}}         % Dinh and Sibony
\newcommand{\FV}{{FV}}         % Falqui and Viallet
\newcommand{\FJ}{{FJ}}         % Favre and Jonsson
\newcommand{\Fa}{{Fa}}         % Favre
\newcommand{\FoRe}{{FR}}       % Fomin and Reading
\newcommand{\FZ}{{FZ}}         % Fomin and Zelevinsky
\newcommand{\Fra}{{Fr1}}       % Friedland 1990 
\newcommand{\Frb}{{Fr2}}       % Friedland 1991
\newcommand{\Frc}{{Fr3}}       % Friedland 1995
\newcommand{\Frd}{{Fr4}}       % Friedland 2006
\newcommand{\FM}{{FM}}         % Friedland and Milnor
\newcommand{\Ga}{{Ga}}         % Gale
\newcommand{\GRP}{{GRP}}       % Grammaticos, Ramani, and Papageorgiou
\newcommand{\Gu}{{Gu}}         % Guedj
\newcommand{\Ha}{{Ha}}         % Hartshorne (Zariski, etc.)
\newcommand{\HK}{{HK}}         % Hasselblatt and Katok
\newcommand{\HNP}{{HNP}}       % Hasselblatt, Nitecki and Propp
\newcommand{\HVa}{{HV1}}       % Hientarinta and Viallet 1
\newcommand{\HVb}{{HV2}}       % Hientarinta and Viallet 2
\newcommand{\Hoa}{{H1}}        % Hone
\newcommand{\Hob}{{H2}}        % Hone
\newcommand{\KH}{{KH}}         % Katok and Hasselblatt
\newcommand{\LRGOT}{{LRGOT}}   % Lafortune, Ramani, Grammaticos, Ohta, 
                               % and Tamizhmani,
\newcommand{\LW}{{LW}}         % Lind and Ward
\newcommand{\Ma}{{Ma}}         % Maegawa,
\newcommand{\Mu}{{Mu}}         % Mumford
\newcommand{\MP}{{MP}}         % Musiker and Propp
\newcommand{\Ng}{{Ng}}         % Nguyen
\newcommand{\Os}{{Os}}         % Osin
\newcommand{\RGLO}{{RGLO}}     % Ramani, Grammaticos, Lafortune, and Ohta
\newcommand{\RS}{{RS}}         % Russakovskii and Shiffman
\newcommand{\St}{{St}}         % Stanley
\newcommand{\Ta}{{T1}}         % Takenawa 1
\newcommand{\Tb}{{T2}}         % Takenawa 2
\newcommand{\TEGORS}{{TEGORS}} % Takenawa, Eguchi, Grammaticos, Ohta, 
                               % Ramani,  and Satsuma
\newcommand{\V}{{V}}           % Veselov
\newcommand{\W}{{W}}           % Wielandt
\newcommand{\Ze}{{Ze}}         % Zelevinsky

% Macros
\newcommand{\co}{{\,:\,}}
\newcommand{\Z}{{\mathbb Z}}
\newcommand{\R}{{\mathbb R}}
\newcommand{\C}{{\mathbb C}}
\newcommand{\PP}{{\mathbb P}}
\newcommand{\Ac}{{\mathcal A}}
\newcommand{\Pc}{{\mathcal P}}
\newcommand{\CP}{{\bf CP}}
\newcommand{\Max}{{\rm Max}}
\newcommand{\dfn}{\mathbin{{:}{=}}}
\makeatletter
\global\def\tlabel#1{\label{NUM#1}\let\@d@mmy\@currentlabel\def\@currentlabel{Theorem~\@d@mmy}\label{#1}}
\global\def\llabel#1{\label{NUM#1}\let\@d@mmy\@currentlabel\def\@currentlabel{Lemma~\@d@mmy}\label{#1}}
\global\def\plabel#1{\label{NUM#1}\let\@d@mmy\@currentlabel\def\@currentlabel{Proposition~\@d@mmy}\label{#1}}
\global\def\clabel#1{\label{NUM#1}\let\@d@mmy\@currentlabel\def\@currentlabel{Corollary~\@d@mmy}\label{#1}}
\global\def\jlabel#1{\label{NUM#1}\let\@d@mmy\@currentlabel\def\@currentlabel{Conjecture~\@d@mmy}\label{#1}}
\global\def\cllabel#1{\label{NUM#1}\let\@d@mmy\@currentlabel\def\@currentlabel{Claim~\@d@mmy}\label{#1}}
\global\def\dlabel#1{\label{NUM#1}\let\@d@mmy\@currentlabel\def\@currentlabel{Definition~\@d@mmy}\label{#1}}
\global\def\elabel#1{\label{NUM#1}\let\@d@mmy\@currentlabel\def\@currentlabel{Example~\@d@mmy}\label{#1}}
\global\def\rlabel#1{\label{NUM#1}\let\@d@mmy\@currentlabel\def\@currentlabel{Remark~\@d@mmy}\label{#1}}
\global\def\qlabel#1{\label{NUM#1}\let\@d@mmy\@currentlabel\def\@currentlabel{Question~\@d@mmy}\label{#1}}
\global\def\flabel#1{\label{NUM#1}\let\@d@mmy\@currentlabel\def\@currentlabel{Figure~\@d@mmy}\label{#1}}
\global\def\Clabel#1{\label{NUM#1}\let\@d@mmy\@currentlabel\def\@currentlabel{Chapter~\@d@mmy}\label{#1}}
\global\def\Slabel#1{\label{NUM#1}\let\@d@mmy\@currentlabel\def\@currentlabel{Section~\@d@mmy}\label{#1}}
\global\def\sbslabel#1{\label{NUM#1}\let\@d@mmy\@currentlabel\def\@currentlabel{Subsection~\@d@mmy}\label{#1}}
\global\def\ssslabel#1{\label{NUM#1}\let\@d@mmy\@currentlabel\def\@currentlabel{Subsubsection~\thesection\thesubsection\@d@mmy}\label{#1}}
\global\def\numref#1{\ref{NUM#1}}

\title{Degree-growth of monomial maps}
\author{Boris Hasselblatt}
\address{Department of Mathematics\\
Tufts University\\
Medford, MA 02155\\
USA}
\email{bhasselb@tufts.edu}
\author{James Propp}
\address{Department of Mathematical Sciences \\
University of Massachusetts Lowell\\
Lowell, MA 01854}
\email{James\_Propp@uml.edu.removelast23characters}
\begin{abstract} 
For projectivizations of rational maps 
Bellon and Viallet defined the notion of algebraic entropy 
using the exponential growth rate of the degrees of iterates. 
We want to call this notion to the attention of dynamicists 
by computing algebraic entropy for 
certain rational maps of projective spaces (\ref{THM2}) 
and comparing it with topological entropy (\ref{THM1}).
The particular rational maps we study are monomial maps (\ref{DEFmonom}),
which are closely related to toral endomorphisms.
Theorems\ \ \numref{THM1} and~\numref{THM2} imply 
that the algebraic entropy of a monomial map
is always bounded above by its topological entropy,
and that the inequality is strict
if the defining matrix has more than one eigenvalue outside the unit circle.
Also, Bellon and Viallet conjectured
that the algebraic entropy of every rational map
is the logarithm of an algebraic integer, 
and\ \ \ref{THM2} establishes this for monomial maps.
However, a simple example (the monomial map of \ref{EXPLbadinverse})
shows that a stronger conjecture of Bellon and Viallet is incorrect,
in that the sequence of algebraic degrees
of the iterates of a rational map of projective space
need not satisfy a linear recurrence relation with constant coefficients.
\end{abstract}
\maketitle
\section{Introduction}\Slabel{SecIntro}
\subsection{Algebraic entropy}
In their 1998 paper \cite{BV}, Bellon and Viallet introduced 
the concept of ``algebraic entropy''
for the study of iterates of rational maps,
measuring the rate at which the algebraic degree
of the $N$th iterate of the map grows as a function of $N$.
This natural and appealing notion
(foreshadowed in work of Arnold \cite{Ar}
and paralleled in work by Russakovskii and Shiffman \cite{RS},
albeit with different terminology)
seems to have escaped the attention of
most researchers in ergodic theory and dynamical systems;
to our knowledge, the only articles on this topic
that have appeared in \emph{Ergodic Theory and Dynamical Systems}
thus far are \cite{Ma} and \cite{Gu}.
Hence, a major motivation behind the writing of this article
is a desire to advertise the study of degree-growth
and to encourage readers of this journal
to think about transporting established ideas
from measurable and topological dynamics
into the setting of algebraic geometry.
More specifically, the following conjecture
deserves attention from dynamicists
of an algebraic bent:
\begin{conjecture}[Bellon and Viallet]\jlabel{CONJBellonViallet}
The algebraic entropy of every rational map
is the logarithm of an algebraic integer.
\end{conjecture}
\subsection{Monomial maps and projectivization}
A second purpose in writing this article is to show that a simple class of
rational maps provides insight into fundamental questions about algebraic
entropy.
\begin{definition}\dlabel{DEFmonom}
Every $n$-by-$n$ nonsingular integer matrix $A=(a_{ij})_{i,j=1}^n$
determines a mapping $(x_1,\dots,x_n)\mapsto(y_1,\dots,y_n)$ from a dense
open subset $U$ of complex $n$-space $\C^n$ to itself by $$y_i = \prod_j
x_j^{a_{ij}}.$$ (If all $a_{ij}\ge0$, then $U=\C^n$.) We call this an
\emph{affine monomial map}.
\end{definition}
\begin{remark}%\rlabel{REM}
The map carries the $n$-torus
$\{(x_1,\dots,x_n): |x_1|=\dots=|x_n|=1\}$ to itself,
and the restriction of the map to the $n$-torus
is isomorphic to the toral endomorphism associated with $A$.
\end{remark}
In this article we focus on
a slightly different construction,
namely, the projectivization of the affine monomial map.
Each projectivized monomial map sends a certain dense open subset $U$
of complex projective $n$-space $\CP^n$ to itself.
(See Section~\SecDef~for relevant definitions and notation.)
Moreover, the action of the map on the $n$-torus 
$\{(x_1 \co \dots \co x_{n+1}): |x_1| = \dots = |x_{n+1}| \neq 0\} 
\subseteq \CP^n$ 
is once again isomorphic to the toral endomorphism associated with $A$.

\begin{remark}%\rlabel{REM}
In accordance with algebraic geometry nomenclature,
we refer to maps from $\C^n$ to itself as ``affine'' 
and maps from $\CP^n$ to itself as ``projective''.
\end{remark}
\begin{example}
\label{soletwo}
Let $A$ be 
the 1-by-1 matrix whose sole entry is 2.
The affine monomial map associated with $A$
is the squaring map $z \mapsto z^2$ on $\C$,
whereas the projective monomial map associated with $A$
is the squaring map on the complex projective line $\CP^1$,
also known as the Riemann sphere $\C \cup \{ \infty \}$.
\end{example}
Example~\ref{soletwo}~is atypical in that
the squaring map is well-defined on all of $\CP^1$.
Later we will see for most integer matrices $A$
we need to restrict the monomial map associated with $A$
to a dense open proper subset $U$ of $\CP^n$.
(From here on, the term ``monomial map'' will usually refer to
a complex projective monomial map unless otherwise specified.)
\subsection{Relations between entropies}
A monomial map restricted to $U$ is continuous,
so it makes sense to ask about its topological entropy.
Since $U$ is typically not a compact space,
it is not immediately clear how the topological entropy should be defined;
fortunately, \cite{HNP} shows that some of the most natural candidate
definitions agree and clarifies the relation between the main notions
that have been proposed.
In \ref{SecTop}, we show that
for this notion of topological entropy, 
the topological entropy of the monomial map
associated with the matrix $A$
is no less than
the topological entropy of the toral endomorphism
associated with $A$,
which in turn is equal to the logarithm of the product of $|z|$
as $z$ ranges over all the eigenvalues of $A$
outside the unit circle
(\ref{THM1}).

At the same time,
monomial maps fall into the framework of Bellon and Viallet,
and we show (\ref{THM2})
that the algebraic entropy of a monomial map
is equal to the logarithm of the spectral radius
of the associated $n$-by-$n$ integer matrix,
i.e., the maximum value of the logarithm of $|z|$
as $z$ ranges over all the eigenvalues of $A$.

Theorems \numref{THM1} and \numref{THM2} 
imply that the algebraic entropy of a monomial map does
not exceed its topological entropy, and that the inequality is strict if
the defining matrix has more than one eigenvalue outside the unit circle.

Since the entries of $A$ are integers, the eigenvalues of $A$ are all
algebraic integers. Thus \ref{THM2} (or, rather, \ref{CORBVmonomial})
provides support for the  Bellon--Viallet \ref{CONJBellonViallet}. On the
other hand, we devise a monomial map that falsifies a stronger conjecture
of Bellon and Viallet's, namely, that the sequence of degrees of the
iterates of a rational map satisfy a linear recurrence with constant
coefficients. The trick is to choose a matrix $A$ whose dominant
eigenvalues are a pair of complex numbers $re^{i \theta}$, $re^{-i \theta}$
where $\theta$ is incommensurable with $2\pi$. For such an $A$, the
sequence of degrees is a patchwork of a finite collection of integer
sequences that individually satisfy linear recurrences with constant
coefficients; the degree sequence jumps around between elements of the
family in a nonperiodic fashion. Details are given in \ref{SecCounter}.

We also describe in Section~\SecPiece\ 
an analogue of algebraic entropy
applicable to the dynamics of piecewise linear maps.

These discoveries are not deep;
they illustrate that there is a lot of ``low-hanging fruit''
in the study of iteration of rational maps
from a projective space to itself,
and suggest that a more vibrant interaction
between the dynamical systems community
and the integrable systems community
(perhaps mediated by researchers in the
field of several complex variables)
could lead to more rapid progress
in the development of the theory
of algebraic dynamical systems.

\section{Definitions}\Slabel{SecDef}
We review some basic facts about projective geometry
(more details can be found in \cite{Mu})
before commencing a discussion of algebraic degree
and algebraic entropy (drawing heavily on \cite{BV}).
\subsection{Projective space}
\begin{definition}
Complex projective $n$-space
is defined as $\CP^n = (\C^{n+1} \setminus \{0\} ) / \sim$,
where $u \sim v$ iff $v = cu$ for some $c\in\C\setminus\{0\}$.
We write the equivalence class of $(x_1,x_2,\dots,x_{n+1})$ 
in $\CP^n$ as $(x_1 \co x_2 \co \dots \co x_{n+1})$.

The standard embedding
$\Pc\colon(x_1,x_2,\dots,x_n) \mapsto (x_1 \co x_2 \co \dots \co x_n \co 1)$ 
of affine $n$-space into projective $n$-space has an ``inverse map''
$\Ac\colon(x_1 \co \dots \co x_n \co x_{n+1})
\mapsto 
(\frac{\textstyle x_1}{\textstyle x_{n+1}},\dots,
\frac{\textstyle x_n}{\textstyle x_{n+1}})$.
The ratios $x_i/x_{n+1}$ ($1 \leq i \leq n$),
defined on a dense open subset of $\CP^n$,
are the \emph{affine coordinate variables} on $\CP^n$.
\end{definition}
Geometrically, one may model $\CP^n$ as the set of lines
through the origin in $(n+1)$-space.  In this model,
the point $(a_1 \co a_2 \co \dots \co a_{n+1})$ in $\CP^n$
corresponds to the line
$x_1 / a_1 = x_2 / a_2 = \cdots = x_{n+1} / a_{n+1}$ in $\C^{n+1}$
(for those $i$ with $a_i = 0$, we impose the condition $x_i = 0$).
The intersection of this line with the hyperplane $x_{n+1}=1$
is the point
$$(\frac{a_1}{a_{n+1}},\frac{a_2}{a_{n+1}},\dots,\frac{a_n}{a_{n+1}},1)$$
(as long as $a_{n+1} \neq 0$).
We identify affine $n$-space with the hyperplane $x_{n+1}=1$.
Affine $n$-space in this way becomes a
Zariski-dense subset of projective $n$-space.
(See e.g., \cite{Ha} for the definition and basic properties
of the Zariski topology.)
Since there is nothing special about the $n+1$st coordinate in $\CP^n$,
each of the hyperplanes $x_i=1$ ($1 \leq i \leq n+1$)
is a copy of (complex) affine $n$-space.
Thus we might see projective $n$-space as the result
of gluing together $n+1$ affine $n$-spaces in a particular way.
Under this viewpoint,
a monomial map is the result of gluing together
$n+1$ compatible toral endomorphisms in a particular way.
\begin{definition}
We define the distance between two points in $\CP^n$
as the angle $0 \leq \theta \leq \pi/2$
between the lines in $\C^{n+1}$ associated with those points;
this gives a metric on $\CP^n$,
and the resulting metric topology 
coincides with the quotient topology on
$(\C^{n+1} \setminus \{0\} ) / \sim$.
\end{definition}
\begin{remark}
There is a more natural distance on projective space,
namely the distance induced by the Riemannian  ``Fubini--Study metric'',
and it may play a role in the analysis of
the topological entropy of monomial maps;
however, we will not pursue this topic here.
\end{remark}
\subsection{Rational maps and projectivization}
We will use the term \emph{rational map} in two different ways:
both to refer to a function from (a Zariski-dense subset of) $\C^n$ 
to $\C^m$ given by $m$ rational functions of the affine coordinate variables,
and to refer to the associated function 
from a Zariski-dense subset of $\CP^n$ to $\CP^m$.
(Henceforth, we will refer to rational maps
``from $\C^n$ to $\C^m$'' or ``from $\CP^n$ to $\CP^m$'',
even though the map may be
undefined on a proper subvariety of the domain.)
That is, a ``rational map'' may be affine or projective,
according to context.
\begin{definition}%\dlabel{DEF}
The \emph{projectivization} of an affine map $f$ is the map $\Pc\circ
f\circ\Ac$, written with cleared fractions,
where $\Pc$ and $\Ac$ are as in
Definition 2.1.
\end{definition}
\begin{example}
The partial function $f\colon x \mapsto 1/x$ on affine 1-space
(undefined at $x=0$) is associated with 
the function $g\colon (x \co y) \mapsto (y \co x)$ on projective 1-space
(defined everywhere).
$f$ is its own inverse on its domain,
while $g$ is its own inverse globally.
\end{example}
The maps in this and the next example (some of them partial functions from
affine $n$-space to itself and some of them partial functions from
projective $n$-space to itself) will be called rational maps, and the
context should make it clear whether we are in the affine setting or the
projective setting. In both settings, we identify functions that agree on a
Zariski-dense set. Under this identification, projectivization commutes
with composition, so, in particular, the $N$th power of the
projectivization of an affine map is identified with the projectivization
of the $N$th power of the map.
\begin{example}
The partial function $f\colon (x,y) \mapsto (1/x,1/y)$ on affine 2-space
(undefined at $xy=0$) is associated with the function $g\colon (x \co y \co
z) \mapsto (yz \co xz \co xy)$ on projective 2-space. $g$ is undefined on
$xy = xz = yz = 0$, and its composition with itself is undefined on the
proper subvariety $xyz = 0$ and is the identity map elsewhere. With the
above identification we can say that $g \circ g$ is the identity map and
say that $g$ is self-inverse.
\end{example}
\begin{definition}%\dlabel{DEF}
A \emph{birational (projective) map} 
is a rational map $f$ from $\CP^n$ to $\CP^n$ 
with a rational inverse $g$ (satisfying $f \circ g = g \circ f =$ 
the identity map on a Zariski-dense subset of $\CP^n$).
\end{definition}
\begin{example}\elabel{EXPLyxy}
The affine map $(x,y) \mapsto (y,xy)$
with inverse $(x,y) \mapsto (y/x,x)$
projectivizes as $f\colon(x \co y \co z) \mapsto (yz \co xy \co z^2)$
with inverse $g\colon(x \co y \co z) \mapsto (yz \co x^2 \co xz)$.
(As a check, note that
$f(g(x \co y \co z)) = ((x^2)(xz) \co (yz)(x^2) \co (xz)^2) = 
(x \co y \co z)$.)
\end{example}
\subsection{Degree}
\begin{lemma}\llabel{LEMmapitopoly}
Every rational map from $\CP^n$ to $\CP^m$ can be written in the form 
$(x_1 \co \dots \co x_{n+1}) \mapsto 
(p_1(x_1,\dots,x_{n+1}) \co \dots \co p_{m+1}(x_1,\dots,x_{n+1}))$ 
where the $m+1$ polynomials $p_1,\dots,p_{m+1}$ are homogeneous
polynomials of the same degree having no joint common factor.
\end{lemma}
\begin{proof}
When we apply $\Ac$, we get $n$ ratios of the affine coordinate variables,
with each ratio homogeneous of degree 0.
When we then apply $f$,
we get $m$ rational functions of the affine coordinate variables,
with each rational function homogeneous of degree 0,
and when we apply $\Pc$, we tack on a 1 at the end of the $n$-tuple,
obtaining an $(n+1)$-tuple.
When we clear denominators,
we multiply all $n+1$ of the rational functions of degree 0
by some homogeneous polynomial,
and when we remove common factors,
we divide them by some homogeneous polynomial.
The end result is an $(n+1)$-tuple of
homogeneous polynomials of the same degree,
having no joint common factor
(although any proper subset of the polynomials
may have some factor in common).
\end{proof}
\begin{definition}%\dlabel{DEF}
The common degree of the polynomials in \ref{LEMmapitopoly} is called
the \emph{degree} of the map.
\end{definition}
\begin{example}
The most familiar case is $n=1$:
the rational function $x \mapsto p(x)/q(x)$ 
(where $p$ and $q$ are polynomials with no common factor)
is associated with the projective map
$(x \co y) = (x/y \co 1)\mapsto(p(x/y)/q(x/y) \co 1) =
(p(x/y) \co q(x/y))$.
The rational functions $p(x/y)$ and $q(x/y)$ are homogeneous of degree 0;
to make them polynomials in $x$ and $y$,
we must multiply through by $y^{\max(\deg p, \deg q)}$.
Hence the degree of the mapping is $\max(\deg p, \deg q)$.
\end{example}
\begin{example}\elabel{EXPLdeg2deg3}
A simple example with $n>1$ is given by the projectivization
of the monomial map $(x,y) \mapsto (y,xy)$ of \ref{EXPLyxy}.
The map $f\colon(x \co y \co z) \mapsto (yz \co xy \co z^2)$
is of degree 2, and its square 
$f^2=f\circ f\colon(x\co y\co z)\mapsto((xy)(z^2)\co(yz)(xy)\co(z^2)^2)=
(xyz^2\co xy^2z\co z^4)=(xyz\co xy^2\co z^3)$ is of degree 3.
\end{example}
\begin{example}
More generally, a 2-by-2 nonsingular integer matrix
$$A = \left( \begin{array}{cc} a & b \\ c & d \end{array} \right)$$
is associated with the affine map
$(x,y)\mapsto(x^a y^b, x^c y^d)$
and with the projective map
$(x \co y \co z) = (x/z \co y/z \co 1)\mapsto((x/z)^a (y/z)^b \co (x/z)^c (y/z)^d \co 1)$.
To make all three entries monomials in $x$, $y$, and $z$,
we multiply them by
$x^{\max(-a,-c,0)}$,
$y^{\max(-b,-d,0)}$ and
$z^{\max(a+b,c+d,0)}$,
so the degree of the mapping is 
$\max(-a,-c,0) + \max(-b,-d,0) + \max(a+b,c+d,0)$.
Applying this to the matrices
$$\begin{pmatrix} 0 & 1 \\ 1 & 1 \end{pmatrix}
\text{ and }
\begin{pmatrix} 1 & 1 \\ 1 & 2 \end{pmatrix}$$
reproduces the calculations of the preceding example.
\end{example}
More generally still, we have:
\begin{proposition}\plabel{PROPdegreemonom}
If $A$ is an $n$-by-$n$ nonsingular matrix
with integer entries $a_{ij}$,
the degree of the projective map associated with $A$ is equal to
\begin{equation}\label{EQDA}
D(A)\dfn\sum_{j=1}^n\Max_{i=1}^n(-a_{ij})+\Max_{i=1}^n(\sum_{j=1}^na_{ij}),
\end{equation}
where $\Max(\dots)\dfn\max(0,\dots)$.
\qed
\end{proposition}
For each fixed $n$, the function $D(\cdot)$, viewed as
a function on the space of all real $n$-by-$n$ matrices,
is continuous and piecewise linear.
That is, the hyperplanes given by all the equations
$a_{ij} = 0$ ($1 \leq i,j \leq n$),
$a_{ij} = a_{i'j}$ ($1 \leq i,i',j \leq n$),
$\sum_{j=1}^n a_{ij} = 0$ ($1 \leq i \leq n$), and
$\sum_{j=1}^n a_{ij} = \sum_{j=1}^n a_{i'j}$ ($1 \leq i,i' \leq n$)
yield a decomposition of $\R^{n^2}$ into chambers
such that for all $A$ within each closed chamber $C$,
we have $D(A) = L_C (A)$
for some linear map $L_C\colon\R^{n^2} \mapsto \R$.
Indeed, the degree of the monomial map associated with $A$
is precisely $\max_C L_C(A)$,
where $C$ varies over all the chambers.
\begin{example}[Degree and birational conjugacy]\elabel{EXPLdegreduc}
If we conjugate the involution
$(x,y) \mapsto (y,x)$
via the birational involution
$(x,y) \mapsto (x,x^2-y)$,
we get the involution
$(x,y) \mapsto (x^2-y,(x^2-y)^2-x)$.
When we projectivize, we get a map
$(x:y:z) \mapsto (x^2z^2-yz^3:(x^2-yz)^2-xz^3:z^4)$ of degree 4
that is conjugate to the map $(x:y:z) \mapsto (y:x:z)$ of degree 1.
This demonstrates the important point
that the degree of a projective map
is not invariant under birational conjugacy. 
However, as we will see in the next subsection,
the rate at which the degree of a projective map
grows under iteration of the map
\emph{is\/} invariant under birational conjugacy.
\end{example}
\subsection{Algebraic entropy}\Slabel{BVentropy}
Bellon and Viallet's notion of algebraic entropy,
like most notions of entropy,
owes its existence to an
underlying subadditivity/submultiplicativity property:
\begin{equation}\label{LEMsubadd}
\deg(f \circ g) \leq \deg(f) \deg(g)
\end{equation}
for all rational maps $f,g$.
This is an easy consequence 
of \ref{LEMmapitopoly};
strict inequality in the lemma holds 
precisely when the compositions of the polynomials 
have some factors in common.
This is the ``reduction-of-degree'' phenomenon.

A first consequence of this inequality 
(via a standard argument; e.g., Proposition 9.6.4 of \cite{KH})
is that
$(1/N)\log \deg(f^{N})$
converges, that is, algebraic entropy is well-defined:
\begin{definition}%\dlabel{DEF}
$\lim_{N\to\infty}({1}/{N}) \log \deg(f^{N})$ is called the
(Bellon--Viallet) \emph{algebraic entropy} of $f$.
\end{definition}
A second consequence of \eqref{LEMsubadd}, no less important,
is that if $g = \phi^{-1} \circ f \circ \phi$
for some birational $\phi$,
then $f$ and $g$ have the same algebraic entropy.
\begin{proposition}%\plabel{PRP}
Algebraic entropy is invariant under birational conjugacy.
\qed
\end{proposition}
\begin{remark}%\rlabel{REM}
It should be mentioned that another use of the term
``algebraic entropy'' occurs in the dynamical systems literature,
measuring the growth of complexity 
of elements of a finitely generated group 
under iteration of some endomorphism of the group;
see, e.g., Definition 3.1.9 in \cite{KH}
and the recent article \cite{Os}.
There does not appear to be any connection
between these two uses of the phrase.
\end{remark}

\section{Existing literature}\Slabel{SecLit}
Bellon and Viallet's definition 
arose from a large body of work
in the integrable systems community
on the issue of degree-growth;
see, e.g., \cite{FV}, \cite{HVa} and \cite{HVb}.
More recent articles on the topic
coming from this community
include \cite{Be}, \cite{LRGOT} and \cite{RGLO}.
\subsection{Dynamical degrees}
A notion equivalent to Bellon and Viallet's was introduced 
at the same time in independent work by 
Russakovskii and Shiffman \cite{RS},
drawing upon earlier work by Friedland and Milnor \cite{FM}.
Russakovskii and Shiffman's theory
associates various quantities, called dynamical degrees,
with a rational map;
the algebraic entropy
is simply the logarithm of the dynamical degree of order 1.
To give the flavor of this work
(without purporting to define the notions being used),
we state that the $k$th dynamical degree of a rational map $f$ 
from $\CP^n$ to itself
is given by $$\lim_{N \rightarrow \infty} 
\left( \int (f^N)^* (\omega^k) \wedge \omega^{n-k} \right)^{1/N}$$
where $\omega$ denotes a K\"ahler form on $\CP^n$
(a complex (1,1) form).
\subsection{Intersections}
Algebraic entropy has antecedents elsewhere in dynamics.
For, as was pointed out by Bellon and Viallet,
the degree of a map
is equal to the number of intersections
between the image of
a generic line in $\CP^n$
and a generic hyperplane in $\CP^n$.
Thus algebraic entropy measures the growth rate
of the number of intersections
between one submanifold
and the image of another submanifold,
and is therefore related to the 
intersection-complexity research program of Arnold \cite{Ar},
introduced in the early 1990s 
and mostly neglected since then by mathematicians
(though studied by some physicists: see e.g., \cite{BM} and \cite{AABM}).
The intermediate dynamical degrees of Russakovskii and Shiffman
can be given definitions in this framework;
specifically, the $k$th dynamical degree of a rational map
from $\CP^n$ to itself (for any $k$ between 0 and $n$)
is equal to the number of intersections
between the image of
a generic $\CP^k$ in $\CP^n$
and a generic $\CP^{n-k}$ in $\CP^n$.
Taking $k=n$,
we see that the top dynamical degree of a rational map
is precisely its topological degree
(the number of preimages of a generic point).

It is worth remarking that some articles
(such as \cite{BM} and \cite{AABM}),
in keeping with Arnold's terminology,
use the term ``complexity'' of $f$ to refer to 
$\lim_{N\to\infty}\sqrt[N]{\deg(f^N)}$,
so complexity is just another name
for dynamical degree of order 1.

More recent articles on the topic of dynamical degree, intersections,
and algebraic entropy include 
\cite{BK}, \cite{BFJ}, \cite{DF}, \cite{DS}, \cite{FJ}, 
\cite{Ta}, \cite{Tb}, and \cite{TEGORS}.
These articles often employ the language of
several complex variables,
with the apparatus of de Rham currents and cohomology.
See also Friedland's survey \cite{Frd}.

Lastly, we mention Veselov's survey article \cite{V},
which contains a good treatment of the multifarious notion
of integrability.

\section{Examples}\Slabel{SecEx}
In this section we present a collection of examples, making some basic
observations about most of them.  Several of these examples appear
repeatedly in later sections to illustrate salient points at appropriate
times.
\begin{example}
The H\'enon map
$(x,y) \mapsto (1+y-Ax^2, Bx)$
projectivizes as
$(x \co y \co z) \mapsto (z^2+yz-Ax^2 \co Bxz \co z^2)$.
For any nonzero constants $A$ and $B$, the $N$th iterate of
this map has degree $2^N$, so every nondegenerate H\'enon map
has algebraic entropy $\log 2$.
This important example is discussed in detail by Bellon and Viallet.
\end{example}
\begin{example}\elabel{EXPLMusiker}
The map $f\colon(x,y) \mapsto (y,(y^2+1)/x)$ is the composition of the two
involutions $(x,y) \mapsto ((y^2+1)/x,y)$ and $(x,y) \mapsto (y,x)$ but is
itself of infinite order. Its projectivization is the map $(x \co y \co z)
\mapsto (xy \co y^2+z^2 \co xz)$. It can be shown that the degree of $f^N$
is only $2N$. Hence the algebraic entropy of $f$ is zero.  This example is 
discussed in greater depth in \cite{MP}, \cite{Ze}, and \cite{Hob}.  
(Amusingly, if one replaces $y^2$ by $y$ in the definition of the affine 
map $f$, one obtains a map of order 5 that was probably known to Gauss 
because of its connection with his \emph{pentagramma mirificum} and is 
described in some detail in \cite{FoRe}.)
\end{example}
\begin{example}[Somos-4 recurrence]\elabel{EXPLSomos}
The map
$(w,x,y,z)\mapsto(x,y,z,(xz+y^2)/w)$
has a similar flavor.
Its $N$th iterate
has degree that grows like $N^2$,
so it too has algebraic entropy zero.
This is the Somos-4 recurrence,
introduced by Michael Somos 
and first described in print by David Gale \cite{Ga}.
\end{example}
\begin{remark}[Laurent phenomenon]%\rlabel{REM}
In the two preceding examples,
the iterates of the map are all Laurent polynomials
(rational functions that can be written as
a polynomial divided by a monomial)
thanks to ``fortuitous'' cancellations that occur
every time one performs a division that
a priori might be expected to yield
a denominator with more than one term.
(For~\ref{EXPLMusiker}, a proof of ``Laurentness''
can be found in~\cite{SZ};
for~\ref{EXPLSomos}, see \cite[Theorem 1.8]{FZ}.)
Fomin and Zelevinsky call this the ``Laurent phenomenon''.
For instance, in the case of example~\ref{EXPLMusiker},
the iterates of the (affine) map involve rational functions of $x$ and $y$
with denominators $x$, $x^2 y$, $x^3 y^2$, $x^4 y^3$, etc.,
even though a priori one would expect denominators
with two or more terms to arise.
Specifically, $(x,y)$ gets mapped to $(y,(y^2+1)/x)$,
which gets mapped to $((y^2+1)/x,(y^4+2y^2+1+x^2)/x^2y$,
and so on.
Indeed, when one iterates $f$ complicated denominators do arise,
but they always disappear when one cancels common factors
between numerators and denominators.
E.g., when one squares $(y^4+2y^2+1+x^2)/x^2y$, adds 1,
and divides by $(y^2+1)/x$,
one expects to see a factor of $y^2+1$ in the denominator,
but the numerator turns out to contain a factor of $y^2+1$ as well,
so that the end result simplifies to the Laurent polynomial
$(y^6+3y^4+3y^2+2x^2y^2+x^4+2x^2+1)/x^3y^2$.
In the projective context, this simplification turns into an instance 
of the reduction-of-degree phenomenon alluded to in \ref{BVentropy}.
Thus the Laurent phenomenon can be seen as an important case
of the reduction-of-degree phenomenon,
where reduction-of-degree applies in a significant way
to all the iterates of the map. 
The Laurent phenomenon has strong connections
to the confinement-of-singularities phenomenon
(see e.g.~\cite{GRP},~\cite{HVa},~\cite{HVb},~\cite{LRGOT}, and~\cite{Ta}).
\end{remark}
\begin{example}
The map $f\colon(w,x,y,z)\mapsto(x,y,z,z(wz-xy)/(wy-x^2))$
does not quite fall under 
the heading of the Laurent phenomenon,
but comes close.
In the iterates of this map,
the denominators are always 
a power of $xz-y^2$ times a power of $wy-x^2$.
The degrees of these iterates are
3, 5, 9, 13, 17, 23, 29, 37, 45, 53, 63, 73, 85, 97, \dots.
This unfamiliar-looking sequence
is actually five quadratic sequences patched together:
$\deg(f^N)=(2/5) N^2 + (6/5) N + c_N$,
where the $c_N$
depend only on the residue class of $N$ modulo 5.
(Indeed, $\deg(f^N)=\lfloor(2N^2 + 6N + 9)/5\rfloor$;
this formula was guessed by us
and proved by A.\ Hone in private correspondence.)
Once again, the algebraic entropy is zero.
\end{example}
\begin{example}[The Scott map]\elabel{EXPLGaleScott}
The map
$f\colon(x,y,z)\mapsto(y,z,(y^2+z^2)/x)$ is
attributed by David Gale \cite{Ga} to Dana Scott.
This too has the Laurent property
(\cite[Theorem 1.10]{FZ})
and it can be shown (see \cite{Hoa}) that
$\deg(f^N)=2, 4, 8, 14, 24, 40, 66, 108, \dots=2(F_N-1)$,
where $F_N$ denotes the Fibonacci numbers.
Hence this map has algebraic entropy 
$\log \frac{1+\sqrt{5}}{2}$.
\end{example}
\begin{example}[Eigentorus]
A different instance of positive entropy,
much closer to the concerns of this article,
is the monomial map $(x,y) \mapsto (y,xy)$ of \ref{EXPLyxy}
associated with the 2-by-2 matrix
$$A =\begin{pmatrix} 0 & 1 \\ 1 & 1 \end{pmatrix}.$$
Here, 
$$
(x,y)\mapsto(y,xy)\mapsto(xy,xy^2)\mapsto(xy^2,x^2y^3)\mapsto(x^2y^3,x^3y^5)
\mapsto(x^3y^5,x^5y^8)\mapsto\dots
$$
The exponents are Fibonacci numbers,
and the map has algebraic entropy $\frac{1+\sqrt{5}}{2}$.

The associated projective map $(x \co y \co z) \mapsto (yz \co xy \co z^2)$
has an ``eigentorus''
$\{(x \co y \co z): |x| = |y| = |z| \neq 0\}$.
One way to think about this eigentorus
is to consider the matrix
$$A' =\begin{pmatrix} 0 & 1 & 1 \\ 1 & 1 & 0 \\ 0 & 0 & 2\end{pmatrix}$$
obtained from $A$ by adjoining 
a column of nonnegative integers at the right,
in such a fashion that all the row-sums are equal to 2.
Let $V$ and $V'$ denote $\C^2$ and $\C^3$, respectively,
and give them their standard bases,
so that $A\colon V\mapsto V$
and $A'\colon V'\mapsto V'$.
The matrix $A'$ has 
$w = (1,1,1)^T$
as an eigenvector, and we mod out by the eigenspace $W$;
the action of $A'$ on the quotient space $V'/W$ is
isomorphic to the action of $A$ on the original 2-dimensional space $V$.
If we now mod out $V$ by the module generated by 
the two standard unit vectors in $\C^2$
(note: not to be confused with modding out $V$
by the subspace the two vectors span!),
corresponding to the fact that
$e^{2 \pi i m + 2 \pi i n} = 1 = e^0$
for all integers $m,n$,
we get a torus on which $A$ acts as an endomorphism.
The same is true in $\C^3$:
additively modding out by multiples of $w = (1,1,1)^T$
corresponds to the projective identification 
(complex dilation) $\sim$ in $\C^3$.
\end{example}
This situation is quite general: 
\begin{proposition}\plabel{PRPisoeigentorus}
For any nonsingular matrix $A$,
the action of the monomial map associated with $A$,
restricted to the eigentorus,
is isomorphic to the toral endomorphism associated with $A$.
\end{proposition}
\begin{proof}
Recall that every monomial map from $\CP^n$ to itself
can be written in the form
$(x_1 \co \dots \co x_{n+1}) \mapsto 
(p_1(x_1,\dots,x_{n+1}) \co \dots \co p_{m+1}(x_1,\dots,x_{n+1}))$
where the $m+1$ polynomials
$p_1,\dots,p_{m+1}$
are homogeneous monomials of the same degree (call it $d$)
having no joint common factor.
We use the exponents of the $n+1$ variables in the $n+1$ monomials
to form an $(n+1)$-by-$(n+1)$ matrix $A'$,
and argue as above.
\end{proof}
There is a subtle but important point here,
namely, that a monomial map may not be well-defined on all of $\CP^n$,
and that even where the monomial map is well-defined,
iterates of the map may not be. 
A brutal way to deal with the problem
is to restrict the monomial map 
to the subset of $\CP^n$ in which
all $n+1$ affine coordinate variables are nonzero.
A more refined way is to restrict attention
to the set $U\dfn\bigcap_{N \geq 1} {\mbox{dom}}(f^{N})$,
the intersection of the domains of 
the iterated maps $f=f^{1},f^{2},f^{3},...$.

\begin{example}
The projective map 
$f\colon (x \co y \co z) \mapsto (yz \co xy \co z^2)$ from \ref{EXPLyxy} is not well-defined at $(1 \co 0 \co 0)$ or $(0 \co 1 \co 0)$,
and the square of this map is not well-defined at $(1 \co 1 \co 0)$.
We could restrict $f$ to the set $\{(x \co y \co z): xyz \neq 0\}$,
since this restricted map is continuous
(and indeed is a homeomorphism),
but we could also restrict to the more inclusive set
$U = \{(x \co y \co z): z \neq 0\}$.
\end{example}
The only truly well-behaved monomial maps
are those for which the matrix $A$
is a positive multiple of some permutation matrix.
In all other cases,
the projective monomial map has singularities:
\begin{example}
Although the affine map 
$(x,y) \mapsto (x,y^2)$
seems to be nonsingular,
it ``really'' has a singularity at infinity,
as we can see when we projectivize it to
$(x \co y \co z) \mapsto (xz \co y^2 \co z^2)$,
which is undefined at $(1 \co 0 \co 0)$.
\end{example}
The typical monomial map has essential singularities;
there is no way to extend the partial function
to a continuous function defined on all of $\CP^n$.

In this respect, projective monomial maps are somewhat reminiscent
of return maps for nonsmooth billiards,
which share the property of being undefined
on a small portion of the space
(corresponding to trajectories
in which the ball goes into a corner).

However, unlike the billiards case,
in which a seemingly innocuous orbit
can be well-defined for millions of steps
and then suddenly hit a corner,
projective monomial maps have fairly tame sets of singularities,
topologically speaking:
\begin{proposition}\plabel{PRPevtlwelldef}
If $f$ is a monomial map from $\CP^n$ to itself,
and $x$ is a point in $\CP^n$ 
for which $x$, $f(x)$, $f^{2}(x)$, \dots, $f^{n}(x)$
are all well-defined, then $f^{N}(x)$ is well-defined
for all $N>2^{n+1}$.
\end{proposition}
\begin{proof}
Each point in $\CP^n$ can be represented by an $(n+1)$-tuple
of 0's and 1's,
where a 1 stands for any nonzero complex number.
Call this the signature of the point.
It is easy to see that the signature of a point 
determines whether the point is in the domain of $f$,
and in the affirmative case,
determines the signature of the image of the point under $f$.
If $f^k(x)$ is well-defined for all $0 \leq k \leq 2^{n+1}$,
then two of the points $f^k(x)$ must have the same signature,
so that the sequence of signatures has become periodic,
and iteration of the map can be continued indefinitely
without fear of hitting the non-point ``$(0:0:\dots:0)$''.

The bound $2^{n+1}$ can actually be replaced by a much
smaller bound on the order of $n^2$, since the way in 
which the signatures evolve over time correspond
to the way in which the entries of the vectors
$v$, $Mv$, $M^2v$, \dots evolve, where $M$ is a
nonnegative matrix and $v$ is a nonnegative vector;
under this correspondence,
0's correspond to positive entries and 1's correspond to zeroes.
For details on the quadratic bound, see \cite{W}.
\end{proof}
\begin{remark}\rlabel{REMdefU}
\ref{PRPevtlwelldef} shows that the set of $x$ in $\CP^n$ for
which the infinite forward $f$-orbit of $x$ is not well-defined is a union
of proper subspaces that form a (usually nonpure, i.e., mixed-dimension)
complex projective subspace arrangement whose complement $U$ is a dense
open subset of $\CP^n$ and is the natural domain on which to investigate 
the topological dynamics of $f$.
\end{remark}
The dynamics of a monomial map on $U$ can be fairly complicated
combinatorially:
\begin{example}
The monomial map
$f\colon (x \co y \co z) \mapsto (xz \co xy \co z^2)$
associated with the matrix
$$\begin{pmatrix} 1 & 0 \\ 1 & 1 \end{pmatrix}$$
and its iterates are well-defined on most, but not all, 
of the complex projective plane $\CP^2$.
The points that lie on the projective line $z=0$
(excepting the point $(0 \co 1 \co 0)$ itself)
are mapped by $f$ to the point $(0 \co 1 \co 0)$,
which is not in the domain of $f$.
Meanwhile, points on the projective line $y=0$
are fixed points of $f$,
except for the point $(1 \co 0 \co 0)$
(where the projective line $y=0$ meets the projective line $z=0$),
which is not in the domain of $f$.
Also, every point on the projective line $x=0$
is mapped by $f$ to the fixed point $(0 \co 0 \co 1)$.

In the terminology of algebraic geometry,
the 1-dimensional subvariety $x=0$
gets blown down to the 0-dimensional subvariety $x=y=0$,
while the 0-dimensional subvariety $x=z=0$
gets blown up to the 1-dimensional subvariety $x=0$
(to see why the latter assertion is true,
consider how $f$ acts at points near $(0 \co 1 \co 0)$).
\end{example}
For a discussion of iteration of rational maps
that attends to blowing up and blowing down
and its implications for degree-growth,
see \cite{BK}.

\section{Topological entropy}\Slabel{SecTop}
\subsection{Choice of entropy}\sbslabel{SSChEnt}
Recall that \ref{REMdefU} introduced the set $U$
as the set of points $x$ such that $f^{N}(x)$ is defined for all $N \geq 1$.
Since this dense open subset of $\CP^n$
inherits the angle-metric from the compact space $\CP^n$,
we can apply the Bowen--Dinaburg definition of topological entropy
\cite{Bo}, \cite{Di} by way of spanning or separated sets.
But it is desirable to have a more intrinsic way
of thinking about the topological dynamics of $f$.
Friedland's approach in such cases (see \cite{Fra}, \cite{Frb}, and \cite{Frc})
is to compactify the dynamical system
inside a countable product of copies of the original space.
Specifically, one identifies the point $x$
with the orbit $(x,f(x),f^{2}(x),\dots)$ in $(\CP^n)^\infty$,
and takes the closure of the set of all such orbits;
this gives a compact space to which the original
Adler--Konheim--McAndrew definition \cite{AKM} can be applied.
The results of \cite{HNP} show that
these two different ways of defining entropy
coincide in the case of monomial maps.
\subsection{Monomials}
\begin{theorem}\tlabel{THM1}
If $A$ is an $n$-by-$n$ nonsingular integer matrix,
the topological entropy of the monomial map
from $\CP^n$ to itself associated with $A$ (as in \ref{DEFmonom})
is at least the logarithm of the modulus of the product of
all the eigenvalues of $A$ outside the unit circle.
\end{theorem}
\begin{proof}
We use the fact that entropy does not increase when one restricts the dynamical
system to an invariant set.
Hence, by \ref{PRPisoeigentorus}, the topological entropy of the
monomial map on $\CP^n$
is at least the topological entropy of the toral endomorphism 
associated with $A$. 
But the topological entropy of a toral endomorphism is the logarithm of the
modulus of the product of the eigenvalues that lie outside the unit circle
(see \cite{LW} for the history of this result).
\end{proof}
\subsection{Conjectured equality}\sbslabel{SSConjEq}
We believe that the topological entropy of a monomial map,
is exactly equal to the quantity in \ref{THM1},
but we have not found a proof of this.

One way to prove this equality would be 
to make use of the intermediate dynamical degrees
mentioned in \ref{SecLit}.
A theorem of Dinh and Sibony \cite{DS}
says that the topological entropy of a map
is bounded above by the logarithm
of the maximal dynamical degree.
If we order the $n$ eigenvalues of $A$ in such a way that
$|\lambda_1| \geq |\lambda_2| \geq \dots \geq |\lambda_n|$,
then it is natural to conjecture that
the $k$th dynamical degree of a monomial map
is equal to $|\lambda_1 \lambda_2 \cdots \lambda_k|$.
(This conjecture is true for $k=n$:
the product of all the eigenvalues
is equal to the determinant of the matrix $A$,
whose absolute value is the degree
of the associated monomial map.
The conjecture is also true for $k=1$:
this is the content of \ref{THM2} below.)
Note that, as $k$ varies,
the maximum value achieved by $|\lambda_1 \lambda_2 \cdots \lambda_k|$
is equal to the modulus of the product
of those eigenvalues that lie outside the unit circle,
which is known to equal the topological entropy
of the toral endomorphism associated with $A$.
Hence, our conjectural formula
for the dynamical degrees of a monomial map,
in combination with the theorem of Dinh and Sibony,
would imply that the topological entropy of a monomial map
is bounded by the topological entropy of
the associated toral endomorphism.
Since the reverse inequality holds as well
(see \ref{SSChEnt}),
the desired equality would follow.

\section{Algebraic entropy}\Slabel{SecAlg}
Recall the formula \eqref{EQDA} for $D(A)$ that gives the degree of the
monomial map associated with the $n$-by-$n$ matrix $A$.
\begin{remark}%\rlabel{REM}
It is easy to see
that composition of affine monomial maps from $\C^n$ to $\C^n$
is isomorphic to multiplication of $n$-by-$n$ matrices.
So computing the degree of the $N$th iterate of a monomial map
is tantamount to computing $D(A^N)$,
and the algebraic entropy of the monomial map is just
$\lim_{N \rightarrow \infty}(1/N)\log D(A^N)$.
\end{remark}
\begin{theorem}\tlabel{THM2}
If $A$ is an $n$-by-$n$ nonsingular integer matrix,
the algebraic entropy of the monomial map
from $\CP^n$ to itself associated with $A$
is equal to the logarithm of the spectral radius of $A$.
\end{theorem}
\begin{proof}
The entries of $A^N$
are $O(r^N)$,
where $r$ is the spectral radius of $A$, so $D(A^N) = O(r^N)$, 
and the algebraic entropy of the map
is at most the logarithm of the spectral radius.
To prove equality
suppose for the sake of contradiction
that $D(A^N) = O(c^N)$ with $1<c<r$.
Replacing $c$ by a larger constant if necessary,
we get $D(A^N) < c^N$ for all sufficiently large $N$.
Recalling the formula for $D(\cdot)$, we conclude from this
that for large $N$, 
every entry of $A^N$ is greater than $-c^N$
and every row-sum of $A^N$ is less than $c^N$.
That is, we now have
upper bounds on the row-sums of $A^N$ 
and on the negatives of the individual entries of $A^N$;
from these, we can derive 
an upper bound on the entries of $A^N$.
For, since each entry of $A^N$ can be written as
the sum of the entries in its row
minus the $n-1$ entries in that row other than itself,
every entry of $A^N$ is less than $nc^N$.
Hence for every unit vector $u$
(whose components all have modulus less than 1),
each component of $A^N u$ has modulus at most 
$n(n c^N) = n^2 c^N$.
Hence the sum of the squares of the entries of $A^N u$
is at most $n(n^2 c^N)^2$,
so the norm of $A^N u$ is at most $n^{5/2}c^N$.
But when $N$ is large enough,
this estimate contradicts the fact 
(a consequence of the Jordan canonical form theorem)
that $A^N$ has a unit eigenvector $u$
for which the norm of $A^N u$ is $r^N$.
\end{proof}
Theorems \numref{THM1} and \numref{THM2} together imply
\begin{corollary}%\clabel{COR}
For any $A$ with two or more (not necessarily distinct) eigenvalues outside
the unit circle, the algebraic entropy of the monomial map associated with
$A$ is strictly less than the topological entropy of the map.
\qed
\end{corollary}
An immediate consequence of \ref{THM2}
is
\begin{corollary}\clabel{CORBVmonomial}
The algebraic entropy of the monomial map
is equal to the logarithm of an algebraic integer.
\end{corollary}
\begin{proof} 
The spectral radius $r$ of $A$ is an algebraic integer: Let $z$ be a
dominant eigenvalue of $A$, so that $r = |z|$. Since $z$ is an algebraic
integer, so is $\overline{z}$, and hence so is $\sqrt{z \overline{z}} = |z|
= r$.
\end{proof}
This establishes \ref{CONJBellonViallet} for monomial maps.
\section{Counterexamples}\Slabel{SecCounter}
\subsection{Entropy gap}
Another easy consequence of \ref{THM2} is
\begin{corollary}%\clabel{COR}
There exist
monomial maps for which the topological entropy is strictly
greater than the algebraic entropy.
\end{corollary}
\begin{proof}
The affine map $(x,y)\mapsto(x^2,y^3)$ has topological entropy $\log 6$ and
algebraic entropy $\log 3$.
\end{proof}
\subsection{Inverses}
One might be tempted to conjecture
that the algebraic entropy of a birational map
is equal to the algebraic entropy of its inverse
(since most notions of entropy
are preserved by inversion).
Toral automorphisms give an easy way
to see that this fails in general,
because the spectral radius of a matrix
that is invertible over $\Z$
is typically not equal to 
the spectral radius of its inverse:
\begin{example}\elabel{EXPLbadinverse}
Let
$$A\dfn\left( \begin{array}{rrr}
-1 &  1 &  0 \\
-1 &  0 &  1 \\
 1 &  0 &  0 \end{array} \right)$$
with associated monomial map
$f\colon(x,y,z) \mapsto (y/x,z/x,x)$.
The characteristic polynomial of $A$ is $t^3+t^2+t-1$,
whose eigenvalues are approximately
$0.54$ and $-0.77\pm1.12i$.
The spectral radius of $A$ is $\sqrt{(-0.77\dots)^2+(1.12\dots)^2}
\approx 1.36$
and the spectral radius of $A^{-1}$ is 
$1/.54\dots \approx 1.84$,
which is the square of the spectral radius of $A$.
Hence the algebraic entropy of the monomial map $f^{-1}$
is twice the algebraic entropy of $f$. 
\end{example}
\subsection{Degree sequence and linear recurrences}
A more subtle conjecture, due to Bellon and Viallet,
is that
for any rational map $f$,
the sequence $(\deg(f^{N}))_{N=1}^{\infty}$ 
satisfies a linear recurrence with constant 
coefficients and leading coefficient 1.

If this were true, it would certainly imply that the algebraic entropy of a
rational map is always the logarithm of an algebraic integer. However, the
map $f$ in \ref{EXPLbadinverse} gives a counterexample to this
claim:
\begin{proposition}\plabel{PRPnorecurrence}
For the rational map $f$ in \ref{EXPLbadinverse}
the sequence $(\deg(f^{N}))_{N=0}^{\infty}=$
1,\penalty-100\hskip1pt
2,\penalty-100\hskip1pt
3,\penalty-100\hskip1pt
4,\penalty-100\hskip1pt
6,\penalty-100\hskip1pt
9,\penalty-100\hskip1pt
12,\penalty-100\hskip1pt
17,\penalty-100\hskip1pt
25,\penalty-100\hskip1pt
33,\penalty-100\hskip1pt
45,\penalty-100\hskip1pt
65,\penalty-100\hskip1pt
85,\penalty-100\hskip1pt
112,\penalty-100\hskip1pt
159,\penalty-100\hskip1pt
215,\penalty-100\hskip1pt
262,\penalty-100\hskip1pt
365,\penalty-100\hskip1pt
524,\penalty-100\hskip1pt
627,\penalty-100\hskip1pt
833,\dots
of degrees does not satisfy any linear recurrence with constant coefficients.
\end{proposition}
To see what is going on with this example on an intuitive level,
let $d_N$ denote the degree of $f^{N}$,
and consider the sequence
$c_N\dfn-$2,\penalty-100\hskip1pt
{\bf2},\penalty-100\hskip1pt
1,\penalty-100\hskip1pt
$-$5,\penalty-100\hskip1pt
{\bf6},\penalty-100\hskip1pt
0,\penalty-100\hskip1pt
$-$11,\penalty-100\hskip1pt
{\bf17},\penalty-100\hskip1pt
$-$6,\penalty-100\hskip1pt
$-$22,\penalty-100\hskip1pt
{\bf45},\penalty-100\hskip1pt
$-$29,\penalty-100\hskip1pt
$-$38,\penalty-100\hskip1pt
{\bf112},\penalty-100\hskip1pt
$-$103,\penalty-100\hskip1pt
$-$47,\penalty-100\hskip1pt
{\bf262},\penalty-100\hskip1pt
$-$318,\penalty-100\hskip1pt
9,\penalty-100\hskip1pt
571,\penalty-100\hskip1pt
$-$898,\dots,
many of whose entries (shown in boldface)
agree with the corresponding entries of the degree sequence for $f$.
$c_N$
is the sum of the entries in the last row of $A^N$
minus 
the sum of the entries on the principal diagonal of $A^N$.
In terms of the notation introduced following
\ref{PROPdegreemonom},
$c_N=L_C (A^N)$ for a particular chamber $C$.
It appears empirically that the sequence of matrices $A, A^2, A^3, \dots$
visits this chamber $C$ infinitely often,
so that $c_N=d_N$ for infinitely many values of $N$.
Certainly \emph{some} chamber
is visited infinitely often,
so for simplicity we will assume that
this particular chamber gets visited infinitely often.
(The analysis given below does not depend
in any essential way on which chamber $C$ is being discussed.)

\begin{proof}
The sequence of $c_N$'s satisfies the linear recurrence
$c_N = c_{N-3} - c_{N-2} - c_{N-1}$
as a consequence of the Cayley--Hamilton theorem
(note that the characteristic polynomial of this recurrence
coincides with the characteristic polynomial of the matrix $A$), 
so the generating function
$\sum_{N=0}^{\infty} c_N x^N$ is the power series expansion
of a rational function of $x$.
If the sequence $d_N$ satisfied some linear recurrence
with constant coefficients,
then the generating function
$\sum_{N=0}^{\infty} d_N x^N$ would also be the power series expansion
of a rational function of $x$.
It would follow that the generating function
$\sum_{N=0}^{\infty} (d_N - c_N) x^N
= 3x^0 + 0x^1 + 2x^2 + 9x^3 + 0x^4 + \dots$
must also be the power series expansions of a rational function of $x$.
It follows from a standard theorem on such expansions
due in various versions to Skolem, Mahler, and Lech
(see e.g., Exercise 3.a in Chapter 4 of \cite{St})
that the set $S$ consisting of those indices $N$ for which $d_N - c_N = 0$
must be eventually periodic,
that is, there must be some union of (one-sided) arithmetic progressions
whose symmetric difference with $S$ is finite.

To see that this cannot happen,
note that $d_N = c_N$ precisely
when several things are simultaneously true
of the matrix $A^N$:
For $i=1,2,3$ the $j,1$-entry of $A^N$ is nonpositive and does not exceed
any other entry in its column,
%the 1,1 entry of $A^N$ is at most zero
%and is less than or equal to
%both of the other entries in the first column of $A^N$;
%the 2,2 entry of $A^N$ is at most zero
%and is less than or equal to
%both of the other entries in the second column of $A^N$;
%the 3,3 entry of $A^N$ is at most zero
%and is less than or equal to
%both of the other entries in the third column of $A^N$;
and the sum of the entries in the third row of $A^N$
is at least zero and is greater than or equal to
both of the other row-sums of $A^N$.
In particular, if $d_N - c_N$ vanishes 
along some arithmetic progression of values of $N$,
the 1,1 entry of $A^N$ must be nonpositive
along some arithmetic progression of values of $N$.

On the other hand, we can use a basic fact from linear algebra
to express the 1,1 entry of $A^N$ as an algebraic function of $N$.
Recall that the set of solutions of 
a homogeneous linear difference equation
with characteristic polynomial $p(t)$
is spanned by the set of sequences of the form $s_N = N^i r^N$
where $r$ is a root of $p(t)$
and $i$ is some nonnegative integer
strictly smaller than the multiplicity of $r$.
In particular, there is an exact formula for the 1,1 entry of $A^N$
of the form
$c_1 \alpha^N + c_2 \beta^N + c_3 \overline\beta^N$,
where $\alpha$ is the real root of
the characteristic polynomial $t^3+t^2+t-1=0$
and $\beta = r e^{i \theta}$ and $\overline\beta = r e^{-i \theta}$ 
are the complex roots.
Since $c_1 \alpha^N + c_2 \beta^N + c_3 \overline\beta^N$
is real for all $N$,
$c_1$ is real and $c_2$ and $c_3$ are
complex conjugates of one another. 
\begin{lemma}\llabel{LEMincomm}
No power of $\beta$ is real,
i.e., $\theta$ is incommensurable with $2\pi$.
\end{lemma}
\begin{proof}
If there were a positive integer $m$ with $\beta^m$ real,
then $\alpha^m$, $\beta^m$, and $\overline\beta^m$
would be the roots of a cubic with rational coefficients
possessing a double root $\beta^m=\overline\beta^m$;
this would imply that $\alpha^m$ and $\beta^m$ are rational.
But $\alpha^m$, like $\alpha$ itself, is an algebraic integer,
so the only way it can be rational is if it is a rational integer;
and this cannot be, since it is a nonzero real number 
with magnitude strictly between 0 and 1.
\end{proof}
\begin{lemma}\llabel{LEMc12nonzero}
The coefficients $c_2$ and $c_3 = \overline{c_2}$
are nonzero.
\end{lemma}
\begin{proof}
If $c_2$ and $c_3$ vanish, the 1,1 entry of $A^N$ is always $c_1 \alpha^N$. 
Taking two different values of $N$ for which
the 1,1 entry of $A^N$ is an integer,
we find that some power of $\alpha$ is rational, and so then is some power
of $\beta$, contradicting \ref{LEMincomm}.
\end{proof}
\ref{LEMincomm} implies that for values of $N$ 
lying in any fixed arithmetic progression,
the (complex) values taken on by $(\beta/r)^N$ are dense in the unit circle,
and (by \ref{LEMc12nonzero}) the (real) values taken on by 
$c_2 (\beta/r)^N + \overline{c_2} (\overline\beta/r)^N$
are dense in some interval centered at 0.
In particular, for values of $N$ in that arithmetic progression,
$c_1 (\alpha/r)^N + c_2 (\beta/r)^N + \overline{c_2} (\overline\beta/r)^N$
will spend a positive fraction of the time
in a ray of the form $(\epsilon,\infty)$ for some $\epsilon>0$.
This means that the 1,1 entry of $A^N$,
being equal to
$c_1 \alpha^N + c_2 \beta^N + \overline{c_2} \overline\beta^N$,
will be positive for infinitely many values of $N$
(and hence at least one)
in our arithmetic progression.
But this contradicts our choice of the arithmetic progression.

Following back the chain of suppositions,
we see that we must conclude that
the sequence $d_0,d_1,d_2,\dots$
does not satisfy any linear recurrence
with constant coefficients,
and our proof is complete.
\end{proof}
More generally, the same reasoning that is given above
shows 
\begin{proposition}%\plabel{PRP}
Let $A$ be any nonsingular $n$-by-$n$ matrix
whose dominant eigenvalues 
are a pair of complex numbers $re^{i \theta}$, $re^{-i \theta}$
where $\theta$ is incommensurable with $2\pi$.
For iterates of the monomial map associated with $A$,
the degree sequence
does not satisfy any linear recurrence
with constant coefficients.
\qed
\end{proposition}
\begin{example}
The 2-by-2 matrix
$$\left( \begin{array}{rr} 1 & 2 \\ -2 & 1 \end{array} \right)$$ associated
with the (nonbirational) rational map $(x,y) \mapsto (xy^2,y/x^2)$ has
eigenvalues $1 \pm 2i$, and the angle between the lines $y=2x$ and $y=-2x$
is irrational (i.e., incommensurable with $\pi$), so we see that the degree
sequence will not satisfy any linear recurrence with constant coefficients.
(This example is similar to Example 1 of \cite{Fa}.)
\end{example}

\subsection{Conjugation}
We have not studied what happens when one starts with a monomial map
and conjugates it via a nonmonomial birational map,
obtaining (in general) a nonmonomial map.
In particular, it seems conceivable that
a suitable nonmonomial conjugate
of the main counterexample of this paper
might be better behaved,
in the sense that its degree sequence
would satisfy a linear recurrence.

It should be emphasized that
the degree sequence associated with a rational map
is \emph{not} invariant under birational conjugacy.
Conjugating the map $f$
may yield a birational map
with a different degree sequence.
Indeed, we saw in \ref{EXPLdegreduc}
that the very first term of the degree sequence,
namely the degree of the map itself,
may change under birational conjugacy.

\subsection{The price of projectivization}\sbslabel{SSMaillard}
Jean-Marie Maillard, in private communication,
has pointed out that if one works in the affine context
and simply studies iterates of the mapping
$f\colon(x,y,z) \mapsto (y/x,z/x,x)$ in \ref{EXPLbadinverse},
one can express the iterates in closed form:
$f^{N}(x,y,z)$
is a triple of monomials,
each of which can be written in the form
$x^{a_N} y^{b_N} z^{c_N}$
where the sequences 
$a_1,a_2,\dots$\,,
$b_1,b_2,\dots$\,, and
$c_1,c_2,\dots$
\emph{do} satisfy linear recurrence relations
with constant coefficients.
(Since there is no projective cancellation going on here,
this is just a matter of ordinary linear algebra,
in multiplicative disguise.)
Maillard suggests through this example that projectivization,
although conceptually compelling, may come at a price.
In particular, the nonrationality of the degree sequence
for iterates of the associated projective map
might be viewed as a result of our insistence
on working in the projective setting
rather than the affine setting.

Note furthermore that projectivization of the affine monomial map
does not usually remove singularities,
and that projectivization takes a seemingly singularity-free map
like $(x,y) \mapsto (x,y^2)$
and tells us that it actually has a singularity at infinity.
\section{Piecewise linear maps}\Slabel{SecPiece}
Although the main focus of this article has been monomial maps,
a general dynamical theory of birational maps would also treat 
more general maps of the sort considered in \ref{SecEx},
such as the Scott map 
$(x,y,z)\mapsto(y,z,(y^2+z^2)/x)$ in \ref{EXPLGaleScott}.
Just as monomial maps are closely associated with linear
maps from $\R^n$ to itself (which in turn are closely associated
with endomorphisms of the $n$-torus), certain nonmonomial maps 
are associated with piecewise linear maps from $\R^n$ to itself.
\subsection{Subtraction-free maps}
We say a map is subtraction-free
if each component of the map can be written as
a subtraction-free expression in the coordinate variables.
E.g., consider the map $f\colon(x,y) \mapsto (x^2+xy+y^2,x^2-xy+y^2)$.
Since $x^2-xy+y^2 = (x^3+y^3)/(x+y)$, both components of $f(x,y)$
can be written in terms of $x$ and $y$ using only addition,
multiplication, and division.  Hence the mapping is subtraction-free.
This implies that the iterates of $f$ can also be expressed using
only addition, multiplication, and division.  The way in which
this leads us to consider piecewise linear maps is that 
the binary operations $(a,b) \rightarrow \max(a,b)$,
$(a,b) \rightarrow a+b$, and $(a,b) \rightarrow a-b$, 
satisfy many of the same properties as
the binary operations $(x,y) \rightarrow x+y$, $(x,y) \rightarrow xy$, 
and $(x,y) \rightarrow x/y$, respectively
(with the additive identity element 0 in the former setting corresponding 
to the multiplicative identity element 1 in the latter setting).
More specifically, all of the simplifications that occur
when one iterates subtraction-free rational maps
are forced to occur when one iterates the associated
piecewise linear maps.
So, for example, the cancellations that permit
the rational map $(x,y) \mapsto (y,(y+1)/x)$
to be of order 5
force the piecewise linear map
$(a,b) \mapsto (b,\max(b,0)-a)$
to be of order 5 as well.

The operation on subtraction-free expressions
that replaces multiplication by addition,
division by subtraction,
and addition by max, or min,
has attracted a good deal of attention lately;
it is known as ``tropicalization'',
and a good introduction to the topic
is \cite{SS}.

\begin{example}
It is interesting to compare
$(x,y) \mapsto (y,(y^2+1)/x)$ from \ref{EXPLMusiker}
with $(a,b) \mapsto (b,\max(2b,0)-a)$.
Iteration of the former map gives rise to the sequence of rational functions
$x$, $y$, $\dfrac{y^2+1}x$, $\dfrac{y^4+x^2+2y^2+1}{x^2 y}$,
$\dfrac{y^6+x^4+2x^2y^2+3y^4+2x^2+3y^2+1}{x^3 y^2}$, $\dots$,
while iteration of the latter map 
gives rise to the sequence of piecewise linear functions
\begin{alignat*}{7}
&\max(&-1a&\hskip1pt&-0b,&\hskip1pt&-1a&\hskip1pt&+2b,&\hskip1pt&-1a&\hskip1pt&-0b)&,\\
&\max(& 0a&\hskip1pt&-1b,&\hskip1pt&-2a&\hskip1pt&+3b,&\hskip1pt&-2a&\hskip1pt&-1b)&,\\
&\max(& 1a&\hskip1pt&-2b,&\hskip1pt&-3a&\hskip1pt&+4b,&\hskip1pt&-3a&\hskip1pt&-2b)&,\\
&\max(& 2a&\hskip1pt&-3b,&\hskip1pt&-4a&\hskip1pt&+5b,&\hskip1pt&-4a&\hskip1pt&-3b)&,
\end{alignat*}
etc.
(Note that the first of these piecewise linear functions
can be written more simply as $\max(-a,-a+2b)$,
but expressing it in a more redundant fashion
brings out the general pattern.)
\end{example}
\begin{example}\elabel{EXPLFibonacci}
It is even more interesting to consider
the piecewise linear analogue of the Scott map
$(x,y,z)\mapsto(y,z,(y^2+z^2)/x)$ 
from \ref{EXPLGaleScott},
which is $(a,b,c) \mapsto (b,c,\max(2b,2c)-a)$.
Iteration of the latter map
gives rise to the sequence of piecewise linear functions
\begin{alignat*}{13}
&\max(& -1a&\hskip1pt& +2b&\hskip1pt&-0c,&\hskip1pt& -1a&\hskip1pt& -0b&\hskip1pt& +2c,&\hskip1pt&-1a&\hskip1pt& +0b&\hskip1pt&+2c,&\hskip1pt&-1a&\hskip1pt& +2b&\hskip1pt&+0c),&\\
&\max(& -2a&\hskip1pt& +3b&\hskip1pt&-0c,&\hskip1pt& -2a&\hskip1pt& -1b&\hskip1pt& +4c,&\hskip1pt& 0a&\hskip1pt& -1b&\hskip1pt&+2c,&\hskip1pt&-2a&\hskip1pt& +3b&\hskip1pt&-0c),&\\
&\max(& -4a&\hskip1pt& +6b&\hskip1pt&-1c,&\hskip1pt& -4a&\hskip1pt& -2b&\hskip1pt& +7c,&\hskip1pt& 0a&\hskip1pt& -2b&\hskip1pt&+3c,&\hskip1pt&-2a&\hskip1pt& +4b&\hskip1pt& -c),&\\
&\max(& -7a&\hskip1pt&+10b&\hskip1pt&-2c,&\hskip1pt& -7a&\hskip1pt& -4b&\hskip1pt&+12c,&\hskip1pt& 1a&\hskip1pt& -4b&\hskip1pt&+4c,&\hskip1pt&-3a&\hskip1pt& +6b&\hskip1pt&-2c),&\\
&\max(&-12a&\hskip1pt&+17b&\hskip1pt&-4c,&\hskip1pt&-12a&\hskip1pt& -7b&\hskip1pt&+20c,&\hskip1pt& 2a&\hskip1pt& -7b&\hskip1pt&+6c,&\hskip1pt&-4a&\hskip1pt& +9b&\hskip1pt&-4c),&\\
&\max(&-20a&\hskip1pt&+28b&\hskip1pt&-7c,&\hskip1pt&-20a&\hskip1pt&-12b&\hskip1pt&+33c,&\hskip1pt& 4a&\hskip1pt&-12b&\hskip1pt&+9c,&\hskip1pt&-6a&\hskip1pt&+14b&\hskip1pt&-7c),&
\end{alignat*}
etc., in which the coefficients can be expressed
in terms of Fibonacci numbers.
The Lipschitz constants of these maps
grow exponentially,
with asymptotic growth rate given by the golden ratio.
\end{example}
\subsection{Lyapunov growth}
More generally,
when one compares a subtraction-free rational recurrence
with its piecewise linear analogue,
one often finds that the growth rate for the Lipschitz constants 
of iterates of the piecewise linear map
(which one can view as a kind of global Lyapunov exponent)
is equal to 
the growth rate for the degrees of
iterates of the rational map.
In fact, every cancellation that occurs
when one iterates the rational map
also occurs when one iterates the piecewise linear map,
so the algebraic entropy of the former is an upper bound on the
logarithm of the global Lyapunov exponent of the latter.
\begin{remark}%\rlabel{REM}
Purists may note that we are modifying
the usual notion of Lyapunov exponent in several respects.
First, we are re-ordering quantifiers.
Ordinarily one looks at the forward orbit
of a specific point $x$,
and sees how the maps $f^{N}$
expand neighborhoods of $x$
with $N$ going to infinity,
and only after defining this limit
does one let $x$ vary over the space as a whole;
here we are taking individual values of $N$
and for each such $N$ we ask
for the largest expansion that
$f^{N}$ can cause on the whole space.
Another difference is that
our piecewise linear maps are not differentiable,
so we are using Lipschitz constants
as a stand-in for derivatives.
\end{remark}
\subsection{PL maps and PL recurrences}
It may seem that we have wandered a bit
from the main themes of this article, 
but the reader may recall
that piecewise linear maps entered the article fairly early on,
via the formula \eqref{EQDA}.
\begin{example}\elabel{EXPLScott}
The affine Scott map
$(x,y,z) \mapsto (y,z,(y^2+z^2)/x)$ from \ref{EXPLGaleScott} (and
\ref{EXPLFibonacci})
gives rise to a sequence of Laurent polynomials
whose denominators are $x^1 y^0 z^0$,
$x^2 y^1 z^0$,
$x^4 y^2 z^1$,
$x^7 y^4 z^2$,
$x^{12} y^7 z^4$,
$x^{20} y^{12} z^7$, \dots
where the exponent-sequence 0, 1, 2, 4, 7, 12, 20, \dots
is associated with iteration of the piecewise linear map
$(a,b,c) \mapsto (b,c,\max(2b,2c)-a)$ we associated with the Scott map in
\ref{EXPLFibonacci}.
\end{example}

\begin{example}
Consider the affine monomial map
$f\colon(x,y,z) \mapsto (y/x,z/x,x)$ of \ref{EXPLbadinverse}
discussed in \ref{SSMaillard}.
If we write $f^{N}(x,y,z)$ as
$$
(p_N^{1}(x,y,z)/q_N^{1}(x,y,z),
 p_N^{2}(x,y,z)/q_N^{2}(x,y,z),
 p_N^{3}(x,y,z)/q_N^{3}(x,y,z))$$
where (for $1 \leq i \leq 3$) 
$p_N^{i}$ and $q_N^{i}$
are monomials with no common factor,
then we can write each sequence
$p_1^{i}, p_2^{i}, p_3^{i}, \dots$
or
$q_1^{i}, q_2^{i}, q_3^{i}, \dots$
in the form
$x^{a_1} y^{b_1} z^{c_1}$, 
$x^{a_2} y^{b_2} z^{c_2}$, 
$x^{a_3} y^{b_3} z^{c_3}$, \dots
where each of the sequences
$a_1,a_2,a_3,\dots$,
$b_1,b_2,b_3,\dots$, and
$c_1,c_2,c_3,\dots$ 
satisfies a linear recurrence.
Indeed, it is possible that the degree sequence
for iterates of the projective monomial map
$(w \co x \co y \co z) \mapsto (wx \co wy \co wz \co x^2)$
(the projectivization of $f$)
satisfies a \emph{piecewise\/} linear recurrence,
but we have not explored this.
(For a simple example of an integer sequence
that satisfies a piecewise linear recurrence
but does not appear to satisfy any linear recurrence
with constant coefficients, consider the sequence
$1,\allowbreak 1,\allowbreak -1,\allowbreak -1,\allowbreak -3,\allowbreak 1,\allowbreak 3,\allowbreak 9,\allowbreak 7,\allowbreak 3,\allowbreak -11,\allowbreak -11,\allowbreak -17,\allowbreak 11,\allowbreak33,\allowbreak 67,\allowbreak 45,\allowbreak 1, \dots$
satisfying the recurrence
$a_n = \max(a_{n-1},a_{n-2})-2a_{n-3}$.)
\end{example}
\subsection{PL projectivization}
As a final note,
we mention that projectivization
has an analogue in the piecewise linear context,
namely, modding out (additively) by multiples of $(1,1,1)$.
\begin{example}
Consider once again
the map $(a,b,c) \mapsto (b,c,\max(2b,2c)-a)$ from Examples
\numref{EXPLFibonacci} and \numref{EXPLScott}.
It sends $(a',b',c')=(a,b,c)+(d,d,d)$,
to $(b,c,\max(2b,2c)-a)+(d,d,d)$, that is, it commutes with adding
constant multiples of $(1,1,1)$, so
we can consider a quotient action
that acts on equivalence classes of triples,
where two triples are equivalent
if their difference is a multiple of $(1,1,1)$.
\end{example}
This quotient construction applies
whenever our piecewise linear map
is ``homogeneous'',
in the sense that there exists a constant $m$
such that
each component of the piecewise linear map
is a max of linear functions,
all of which have coefficients adding up to $m$.
(In the example we just considered, $m=1$.)

\section{Comments and open questions}\Slabel{SecComments}
We suggest that in some respects, the logarithm of the maximal dynamical
degree behaves in a fashion more analogous with other kinds of entropy 
than Bellon and Viallet's notion of algebraic entropy does.  (Some of our 
e-mail correspondents have taken this point of view as well.) In the case
of a monomial map associated with a nonsingular integer matrix $A$, we have
already shown that algebraic entropy as defined by Bellon and Viallet is
the spectral radius of $A$, whereas the logarithm of the maximal dynamical
degree of the map stands a decent chance of being equal to the topological
entropy of the toral endomorphism associated with $A$. Furthermore,
Tien-Cuong Dinh has pointed out to us in private correspondence that if $f$
is any birational map from projective $n$-space to itself, the $k$th
dynamical degree of $f$ is equal to the $n-k$th dynamical degree of
$f^{-1}$ (as a trivial consequence of the equality between $\int (f^N)^*
(\omega^k) \wedge \omega^{n-k}$ and $\int \omega^k \wedge (f^{-N})^*
\omega^{n-k}$ obtained by a coordinate change), from which it easily
follows that the logarithm of the maximal dynamical degree of $f^{-1}$
equals the logarithm of the maximal dynamical degree of $f$. 
\begin{question}\qlabel{Q1}
Is the algebraic entropy of a monomial map always equal to 
the topological entropy of the associated toral endomorphism?
\end{question}

\begin{question}\qlabel{Q2}
Is the algebraic entropy of a map always bounded above by its topological
entropy?
\end{question}

We have seen that this is true for monomial maps. The discussion in
\ref{SSConjEq} is pertinent.  Also see \cite{Ng}.

A different sort of question about inequalities is:

\begin{question}\qlabel{Q3}
Is algebraic entropy nonincreasing under factor maps?

That is, if we have birational maps 
$f\colon\CP^n \mapsto \CP^n$ and
$g\colon\CP^m \mapsto \CP^m$,
and a rational map $\phi\colon\CP^n \mapsto \CP^m$
satisfying
$$\phi \circ f = g \circ \phi,$$
must the algebraic entropy of $g$
be less than or equal to the algebraic entropy of $f$?
\end{question}
To avoid trivial counterexamples, we should insist
that the map be dominant (i.e., that its image
is Zariski-dense in $\CP^m$);
here, this is equivalent to assuming $n \geq m$.

Of continuing importance is the \ref{CONJBellonViallet} of
Bellon and Viallet:
\begin{question}\qlabel{Q4}
Is the algebraic entropy of a rational map
always the logarithm of an algebraic integer?
\end{question}
One might also try to clarify the situation for the case in which
algebraic entropy vanishes.
\begin{question}\qlabel{Q5}
Can the degree sequence of a rational map
be subexponential but superpolynomial?
\end{question}
\begin{question}\qlabel{Q6}
If the degree sequence of a rational map
is bounded above by a polynomial,
must it grow like $N^k$ for some nonnegative integer $k$,
or can it exhibit intermediate asymptotic behavior, 
such as $\sqrt{N}$?
\end{question}
Even though monomial maps provide counterexamples
to Bellon and Viallet's conjecture about degree sequences,
it surely cannot be a mere coincidence
that so many of the examples studied by Bellon and Viallet and others
have the property that the degree sequences
satisfy recurrence relations with constant coefficients.
So one might inquire whether we can 
rescue Bellon and Viallet's conjecture on degree sequences
by adding extra hypotheses.
One such possible extra hypothesis is suggested by the fact
(pointed out to us by Viallet)
that many of the birational mappings 
studied by Bellon and Viallet
can be written as compositions of involutions.
\begin{question}\qlabel{Q7}
If a rational map is a composition
of involutions, must its degree sequence
satisfy a linear recurrence with constant coefficients?
\end{question}
It may be worth mentioning that, under the hypothesis of \ref{Q7},
the rational map is birationally conjugate to its inverse,
so the two maps have the same algebraic entropy.
\begin{question}\qlabel{Q8}
Must the degree sequence
of a rational map
satisfy a piecewise linear recurrence with constant coefficients?
\end{question}
\begin{question}\qlabel{Q9}
Is there a simple formula
for the intermediate dynamical degrees of monomial maps,
generalizing \ref{PROPdegreemonom}?
\end{question}
Intermediate dynamical degrees (first defined in \cite{RS}),
although conceptually quite natural,
have proved to be difficult to compute in all but the simplest of cases;
monomial maps constitute a setting in which
one might hope to do computations and prove nontrivial results.
It is natural to conjecture that
the $k$th dynamical degree of a monomial map
is equal to $| \lambda_1 \lambda_2 \cdots \lambda_k|$,
where $\lambda_1, \lambda_2, \dots$
are the eigenvalues of the associated matrix,
ordered so that $|\lambda_1| \geq |\lambda_2| \geq \dots$.
As was remarked in \ref{SSConjEq},
a proof of this conjecture for all $k$
would yield an affirmative answer to \ref{Q1}.

\section*{Acknowledgments}
The authors acknowledge generous assistance from
Dan Asimov,
Eric Bedford,
Mike Boyle,
Tien-Cuong Dinh,
Noam Elkies,
Charles Favre,
Shmuel Friedland,
Vincent Guedj,
Andrew Hone,
Kyounghee Kim,
Michael Larsen,
Doug Lind,
Jean-Marie Maillard,
Zbigniew Nitecki,
Nessim Sibony,
Richard Stanley,
Hugh Thomas,
Claude Viallet and the referee, whose comments helped improve the
exposition. 
This work was supported in part
by a grant from the National Security Agency's
Mathematical Sciences Program.

\end{document}